\documentclass[11pt,reqno]{amsart}

\usepackage{cite,a4wide,amsmath,amssymb,graphicx,enumitem,amsaddr,tikz,pgfplots,hhline,epstopdf,multirow,multicol,soul}
\usepackage[hang,flushmargin]{footmisc}

\usepackage{tabu}

\usetikzlibrary{decorations.markings,patterns}

\setlist[enumerate]{leftmargin=8mm}
\setlist[itemize]{leftmargin=8mm}

\newcommand{\cR}{\mathcal{R}}

\newcommand{\ds}{\frac{\text{d}S}{\text{d}t}}
\newcommand{\di}{\frac{\text{d}I}{\text{d}t}}
\newcommand{\dr}{\frac{\text{d}R}{\text{d}t}}

\DeclareMathOperator{\BT}{BT}
\DeclareMathOperator{\GH}{GH}
\DeclareMathOperator{\HB}{HB}
\DeclareMathOperator{\SN}{SN}
\DeclareMathOperator{\PP}{P}

\definecolor{darkgreen}{RGB}{34,177,76}

\begin{document}
\title{Codimension-two bifurcations of an SIR-type model for COVID-19 and their epidemiological implications}

\author{\small Livia Owen, Jonathan Hoseana, Benny Yong}
\address{\normalfont\small Center for Mathematics and Society, Department of Mathematics, Parahyangan Catholic University, Bandung 40141, Indonesia}
\email{benny\_y@unpar.ac.id\textnormal{, }j.hoseana@unpar.ac.id\textnormal{, }livia.owen@unpar.ac.id}
\date{}

\begin{abstract}
We study the codimension-two bifurcations exhibited by a recently-developed SIR-type mathematical model for the spread of COVID-19, as its two main parameters \linebreak ---the susceptible individuals' cautiousness level and the hospitals' bed-occupancy rate--- vary over their domains. We use AUTO to generate the model's bifurcation diagrams near the relevant bifurcation points: two Bogdanov-Takens points and two generalised Hopf points, as well as a number of phase portraits describing the model's orbital behaviours for various pairs of parameter values near each bifurcation point. The analysis shows that, when a backward bifurcation occurs at the basic reproduction threshold, the transition of the model's asymptotic behaviour from endemic to disease-free takes place via an unexpectedly complex sequence of topological changes, involving the births and disappearances of not only equilibria but also limit cycles and homoclinic orbits. Epidemiologically, the analysis confirms the importance of a good control of the values of the aforementioned parameters for a successful eradication of COVID-19. We recommend a number of strategies by which this may be achieved.


\end{abstract}

\maketitle

\section{Introduction}\label{section:Introduction}

The story of COVID-19 is not yet complete. After successfully maintaining an essentially disease-free status for almost two years, China is once again implementing lockdowns, following the unprecedented omicron outbreak, which is mentioned to be ``ten times more severe'' than the original Wuhan outbreak in 2020 \cite{Yong}. Indeed, the country's previously-unchanging maximum number of daily new cases, 15,133, recorded on 13 February 2020, was surpassed on 5 April 2022 with 16,649 new cases, before the latest maximum of 53,345 new cases was reported on 15 April 2022 \cite{JHU}.

The scientific impact of COVID-19 has been remarkable. Over the last three years, the literature has witnessed a surge of interest in the study of the disease's spread, particularly via mathematical models. In mid 2021, we initiated our study by developing the following simple, SIR-type model which incorporates as key parameters the susceptible individuals' cautiousness level $\gamma\in[0,1]$ and the hospitals' bed-occupancy rate $\rho\in[0,1]$:
\begin{equation}\label{eq:model}
\left\{\begin{array}{rcl}
\displaystyle\ds &=& \displaystyle\lambda - \mu S - \frac{\beta SI}{1+\gamma S},\\[0.33cm]
\displaystyle\di &=& \displaystyle-\mu I - \mu' I +\frac{\beta SI}{1+\gamma S}- \frac{\alpha I}{1+\rho I},\\[0.33cm]
\displaystyle\dr &=& \displaystyle-\mu R + \frac{\alpha I}{1+\rho I},
\end{array}\right.
\end{equation}
where $S=S(t)$, $I=I(t)$, and $R=R(t)$ denote the sizes of the susceptible, infected, and recovered subpopulations at time $t\geqslant 0$, while $\beta$, $\lambda$, $\mu$, $\mu'$, and $\alpha$ are positive parameters \cite{YongOwenHoseana}. Subsequently, we applied the model \eqref{eq:model} to the case of Jakarta, with the aim of constructing a quantitative method to determine the appropriate level(s) of social restrictions to be enforced in the region on any given day, based on the latest values of the bed-occupancy rate and the effective reproduction number \cite{YongHoseanaOwen1}. Most recently, as the Indonesian government intensifies its eradicative effort through five forms of interventions: vaccinations, social restrictions, tracings, testings, and treatments, we proposed a substantial modification of the model which takes these into account, with the aim of identifying optimal intervention strategies \cite{YongHoseanaOwen2}.

From the analysis presented in \cite{YongOwenHoseana}, we have seen that the model \eqref{eq:model}, despite its simplicity, exhibits rich dynamical behaviour. Firstly, the model possesses a unique endemic equilibrium $\mathbf{e}_0=(\lambda/\mu,0,0)$ for every set of parameter values, which is stable (unstable) if $\cR_0<1$ ($\cR_0>1$), where
\begin{equation}\label{eq:R0}
\cR_0=\frac{\beta\lambda}{(\mu+\gamma\lambda)\left(\mu+\mu'+\alpha\right)}
\end{equation}
is the model's basic reproduction number, as well as at most three positive endemic equilibria $\mathbf{e}_1$, $\mathbf{e}_2$, $\mathbf{e}_3$. Furthermore, fixing the parameter values
\begin{equation}\label{eq:parameters}
\beta=0.05,\qquad\lambda=10,\qquad\mu=0.01,\qquad\mu'=0.1,\qquad\alpha=0.2,
\end{equation}
$\rho=0.1$ while letting $\gamma$ vary over $[0,1]$, we observed that the model undergoes a number of \textit{codimension-one} bifurcations: backward transcritical, Hopf, and saddle-node
bifurcations of equilibria, as well as homoclinic and saddle-node bifurcations of limit cycles, the latter two being detected via numerical continuation, using AUTO. With regards to the model's \textit{codimension-two} bifurcations, however, we have only pointed out without details in \cite[section 4]{YongOwenHoseana} that, by letting both $\gamma$ and $\rho$ vary over $[0,1]$, one finds instances of Bogdanov-Takens and generalised
Hopf bifurcations. In the present paper, we shall continue the study of the model \eqref{eq:model} by discussing these bifurcations in greater detail, along with their epidemiological implications.


The discussion is organised as follows. In the upcoming section \ref{sec:Overview}, we establish a connection between what has been done in \cite{YongOwenHoseana} and what is to be done in the present paper. We also describe the way in which we use AUTO to detect the aforementioned bifurcations, and give a brief summary of the topological changes occurring near each bifurcation point. In the subsequent section \ref{sec:Details}, we give a more detailed explanation on these changes and what they epidemiologically imply. Essentially, these changes can be viewed as complex manners in which the model's asymptotic behaviour transitions from endemic to disease-free, which involves the births and disappearances of limit cycles and homoclinic orbits, all occurring under the condition that $\cR_0<1$ and that the model's transcritical bifurcation taking place at the basic reproduction threshold is backward. We also recommend several strategies for the disease's eradication which arise from our findings. In the final section \ref{sec:Conclusions}, we summarise our conclusions and describe possible avenues for further investigation.

\section{Overview}\label{sec:Overview}

\begin{figure}\centering
\includegraphics[width=0.7\textwidth]{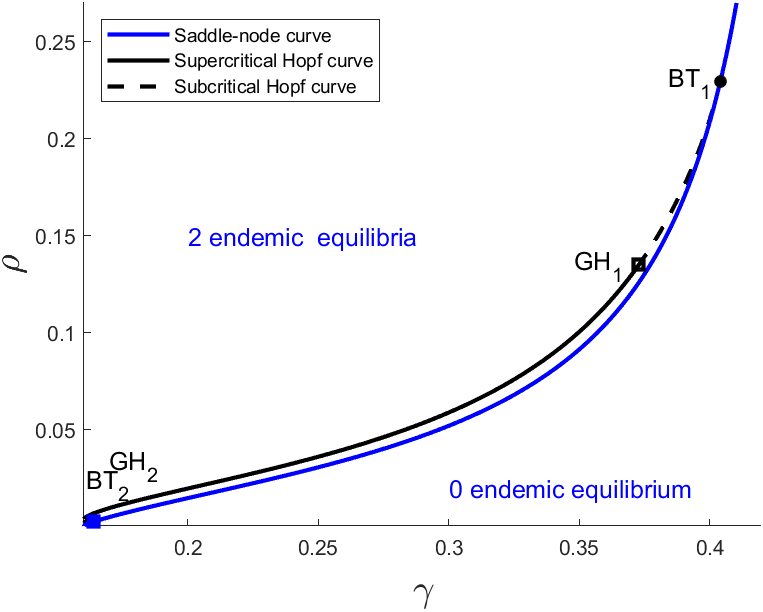}
		\caption{\label{fig:BDglobal}The codimension-two bifurcation diagram of the model \eqref{eq:model} in the region $\left[\gamma_0,0.42\right]\times\left[\rho_0,0.27\right]$ on the $\gamma\rho$-plane, using the values of $\beta$, $\gamma$, $\mu$, $\mu'$, and $\alpha$ listed in \eqref{eq:parameters}.}
\end{figure}

For the rest of the paper, we fix the values of $\beta$, $\lambda$, $\mu$, $\mu'$, and $\alpha$ as in \eqref{eq:parameters}. The basic reproduction number \eqref{eq:R0}, being independent of $\rho$, reduces to a univariate function of $\gamma$:
\begin{equation}\label{eq:R0new}
\cR_0=\frac{5000}{31+31000\gamma}.
\end{equation}
Letting both $\gamma$ and $\rho$ vary over $[0,1]$, we have detected using AUTO \cite{Doedeletal} a set of \textit{bifurcation curves} on the $\gamma\rho$-plane, each of which being a set of points $(\gamma,\rho)$ on the unit square at which the model undergoes a specific bifurcation. In Figure \ref{fig:BDglobal}, we display these curves in the region containing the richest discovered dynamical behaviour: $$\left[\gamma_0,0.42\right]\times\left[\rho_0,0.27\right],\qquad\text{where}\qquad \gamma_0:=\frac{4969}{31000}\quad\text{and}\quad \rho_0:=\frac{29791}{100000000}.$$
As easily verified, in the entire region we have from \eqref{eq:R0} that $\cR_0<1$, and from \cite[Theorem 2.4]{YongOwenHoseana} that the transcritical bifurcation at the basic reproduction threshold is backward. The region, therefore, consists of two adjacent subregions, in each of which the model possesses zero and two endemic equilibria, separated by a \textit{saddle-node bifurcation curve}, containing points $(\gamma,\rho)$ at which these equilibria coalesce. In Figure \ref{fig:BDglobal}, this curve is plotted in blue, and is obtained by carrying out bidirectional continuation beginning from the saddle-node bifurcation point discussed in \cite[section 3]{YongOwenHoseana}: $$\left(\gamma^{(\SN)},0.1\right),\qquad\text{where}\qquad\gamma^{(\SN)}\approx 0.356902.$$ The curve plotted in black, on the other hand, is a \textit{Hopf bifurcation curve}, obtained similarly from the Hopf bifurcation point $$\left(\gamma^{(\HB)},0.1\right),\qquad\text{where}\qquad\gamma^{(\HB)}\approx 0.349638,$$ discussed in \cite[section 3]{YongOwenHoseana}.

The Hopf curve has its endpoints lying on the saddle-node curve:
$$\BT_1\approx\left(0.404023, 0.229494\right)\qquad\text{and}\qquad \BT_2\approx\left(0.164201, 0.002600\right).$$
Letting $\gamma$ and $\rho$ vary smoothly so that the point $(\gamma,\rho)$ travels anticlockwise around each $\BT_i$, one observes the following topological changes, to be detailed in the next section:
\begin{enumerate}
\item[(i)] a homoclinic orbit emanates around a saddle endemic equilibrium via a homoclinic bifurcation, before shrinking and becoming an unstable limit cycle which surrounds a stable endemic equilibrium;
\item[(ii)] the unstable limit cycle is absorbed by the stable endemic equilibrium, which then becomes unstable, via a Hopf bifurcation;
\item[(iii)] the two equilibria coalesce and disappear via a saddle-node bifurcation.
\end{enumerate}
Therefore, at each $\BT_i$, the model undergoes a \textit{Bogdanov-Takens bifurcation} \cite[section 8.4]{Kuznetsov}.

As also apparent in Figure \ref{fig:BDglobal}, the Hopf curve consists of a solid \textit{supercritical} Hopf curve, which indicates the ejection of a \textit{stable} limit cycle, connected at its endpoints
$$\GH_1\approx\left(0.372814, 0.134955\right)\qquad\text{and}\qquad \GH_2\approx\left(0.163907, 0.002496\right)$$
to two \textit{subcritical} Hopf curves, which indicate the ejection of \textit{unstable} limit cycles. As the point $(\gamma,\rho)$ travels anticlockwise around each $\GH_i$, the following topological changes occur, again to be detailed in the next section:
\begin{enumerate}
\item[(i)] a homoclinic orbit emanates around a saddle endemic equilibrium via a homoclinic bifurcation, before shrinking and becoming an unstable limit cycle which surrounds a stable endemic equilibrium;
\item[(ii)] the stable endemic equilibrium loses stability while ejecting a stable limit cycle, via a Hopf bifurcation;
\item[(iii)] the two limit cycles coalesce and disappear, via a saddle-node bifurcation of limit cycles.
\end{enumerate}
Therefore, at each $\GH_i$, the model undergoes a \textit{generalized Hopf bifurcation} \cite[section 8.3]{Kuznetsov}.

We thus have four codimension-two bifurcation points of the model \eqref{eq:model}: $\BT_1$, $\BT_2$, $\GH_1$, and $\GH_2$. In the next section, we shall look at the neighbourhoods of these points, and describe the qualitatively different dynamical behaviours which may be possessed by the model's orbits for various pairs $(\gamma,\rho)$ belonging to these neighbourhoods. From the perspective of epidemiology, these behaviours will confirm the significance of the bifurcation parameters $\gamma$ and $\rho$ for the eradication of COVID-19. The specific epidemiological implications will also be discussed, along with a number of concrete recommendations for the disease's eradication.



\section{Local behaviour and epidemiological implications}\label{sec:Details}

\begin{figure}\centering
{\tabulinesep=0.5mm
\begin{tabu}{|c|c|}\hline
		\includegraphics[width=0.45\textwidth]{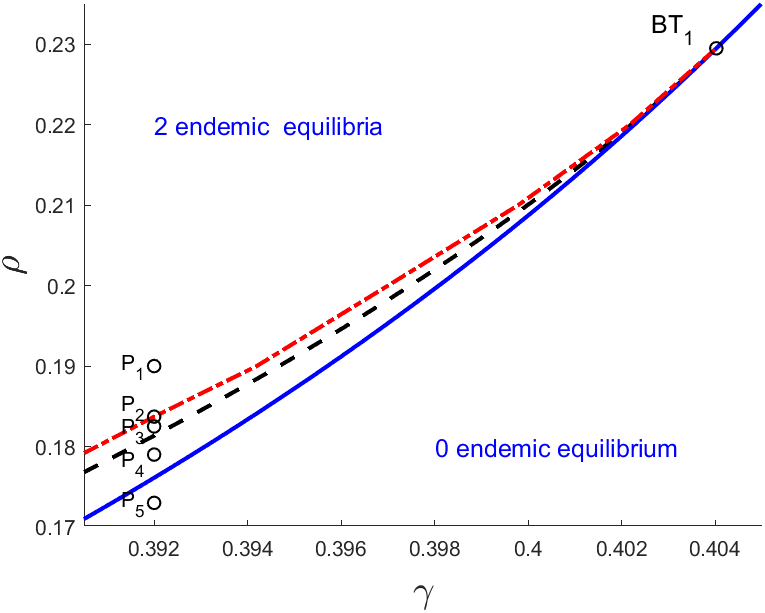}&
		\includegraphics[width=0.45\textwidth]{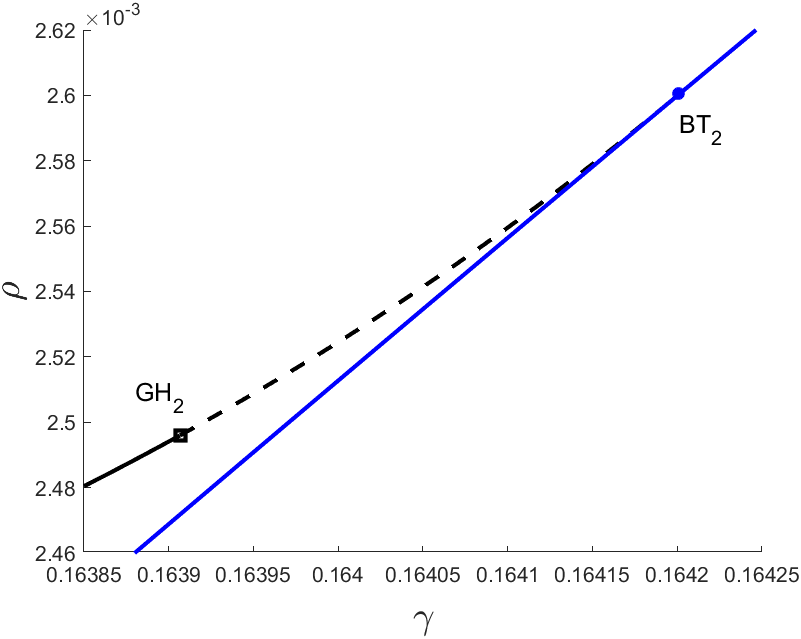}\\
		{\scriptsize \textbf{(a)} Magnification near $\BT_1$} & {\scriptsize \textbf{(b)} Magnification near $\BT_2$}\\\hline
		
		\includegraphics[width=0.45\textwidth]{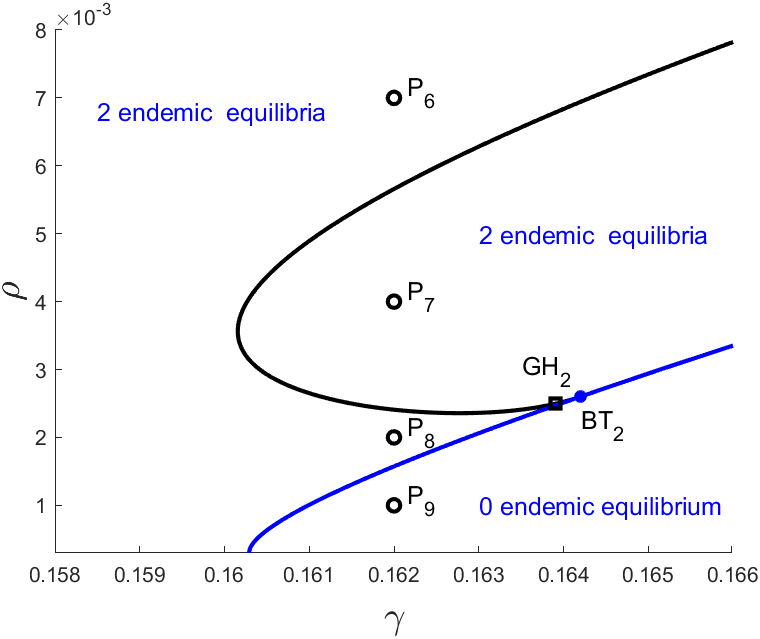}&
		\includegraphics[width=0.45\textwidth]{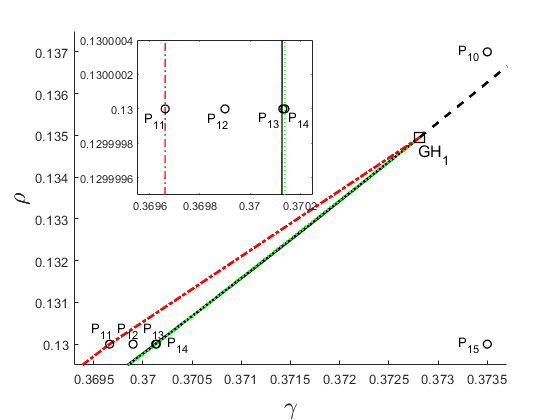}\\
		{\scriptsize \textbf{(c)} Magnification near $\GH_2$} & {\scriptsize \textbf{(d)} Magnification near $\GH_1$}\\\hline
\end{tabu}}\medskip

	\includegraphics[width=0.3\textwidth]{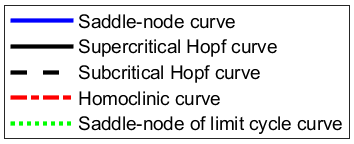}
		\caption{\label{fig:BDlocal}Magnifications of  Figure \ref{fig:BDglobal} near the Bogdanov-Takens points $\BT_1$, $\BT_2$ and generalised Hopf bifurcation points $\GH_1$, $\GH_2$, with additions of homoclinic and saddle-node of limit cycle curves.}
\end{figure}
\begin{figure}\centering
{\tabulinesep=0.5mm
\begin{tabu}{|c|c|c|}\hline
		\includegraphics[width=0.3\textwidth]{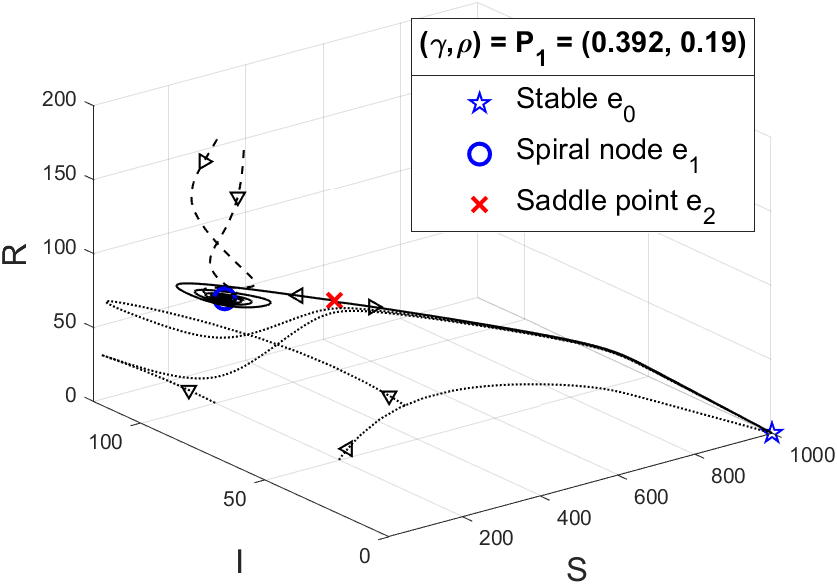}&
		\includegraphics[width=0.3\textwidth]{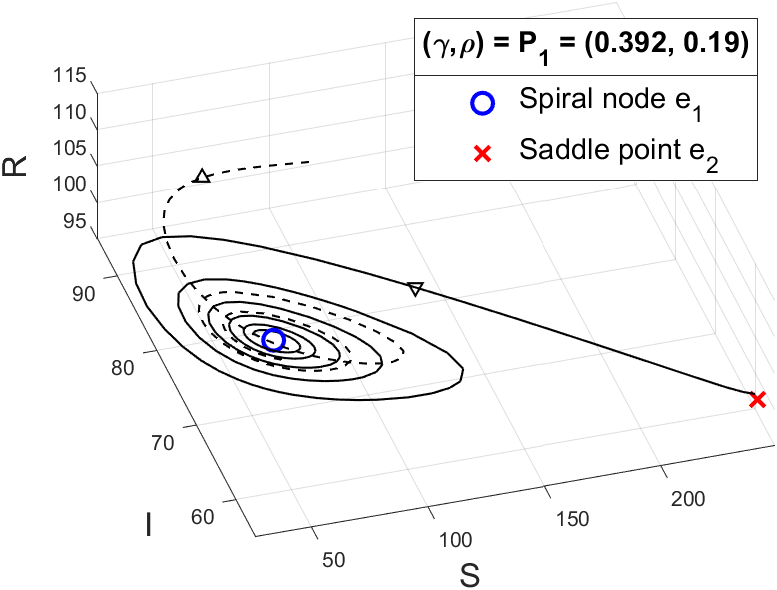}&
		\includegraphics[width=0.3\textwidth]{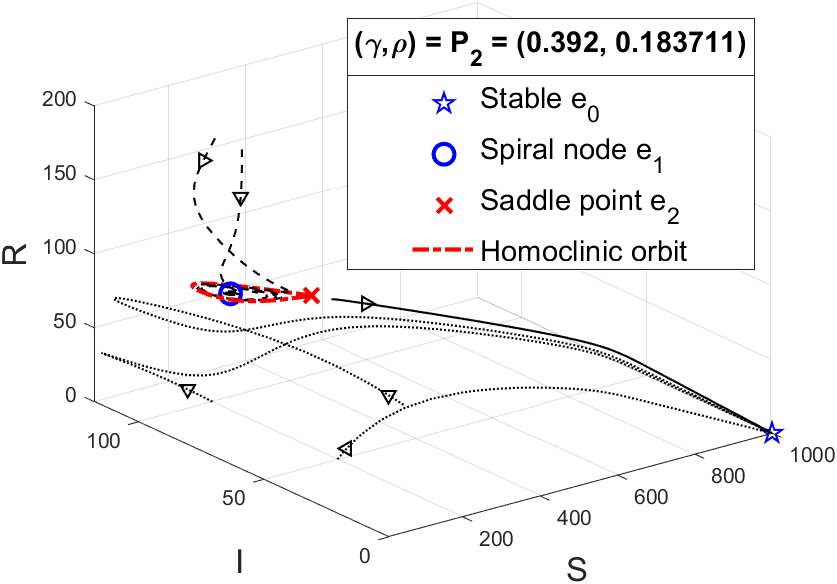}\\
		{\scriptsize \textbf{(a)} $(\gamma,\rho) = \PP_1$} & {\scriptsize \textbf{(b)} $(\gamma,\rho) = \PP_1$ (magnified)} & {\scriptsize \textbf{(c)} $(\gamma,\rho) = \PP_2$}\\\hline
		
		\includegraphics[width=0.3\textwidth]{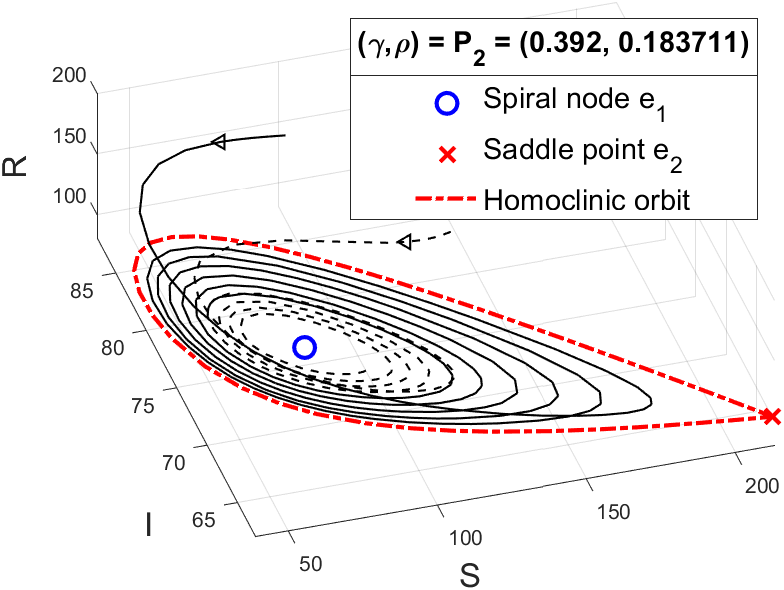}&
		\includegraphics[width=0.3\textwidth]{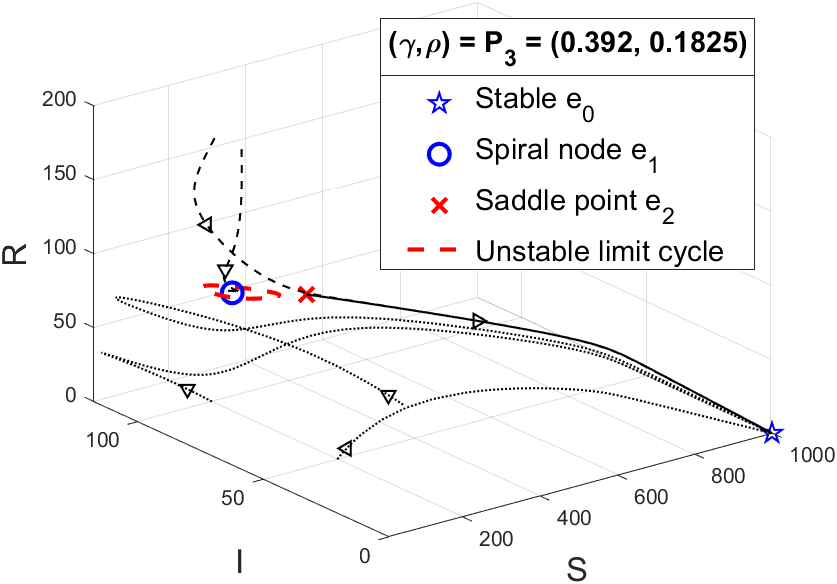}&
		\includegraphics[width=0.3\textwidth]{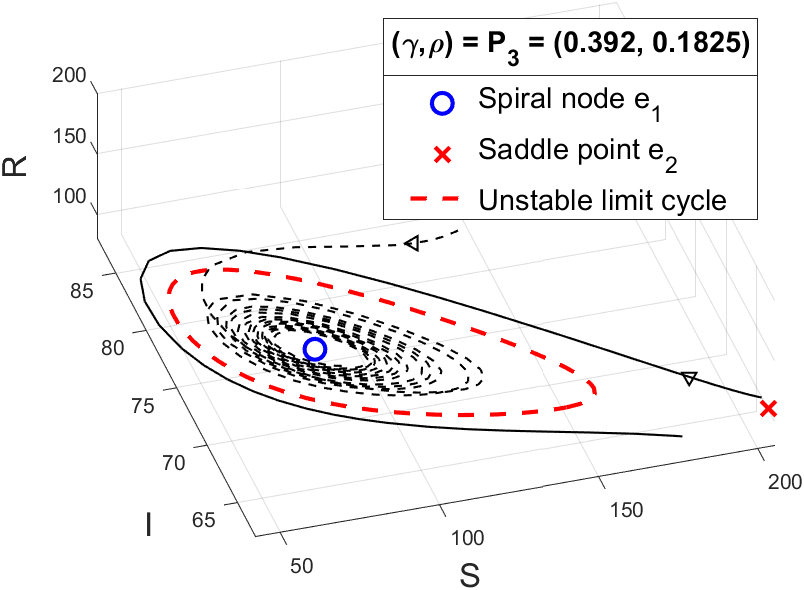}\\
		{\scriptsize \textbf{(d)} $(\gamma,\rho) = \PP_2$ (magnified)} & {\scriptsize \textbf{(e)} $(\gamma,\rho) = \PP_3$} & {\scriptsize \textbf{(f)} $(\gamma,\rho) = \PP_3$ (magnified)}\\\hline

		\includegraphics[width=0.3\textwidth]{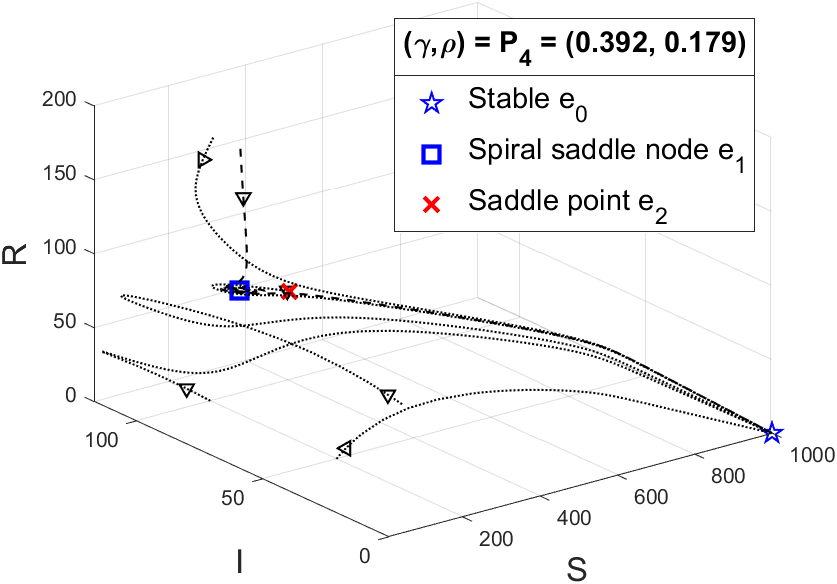}&
		\includegraphics[width=0.3\textwidth]{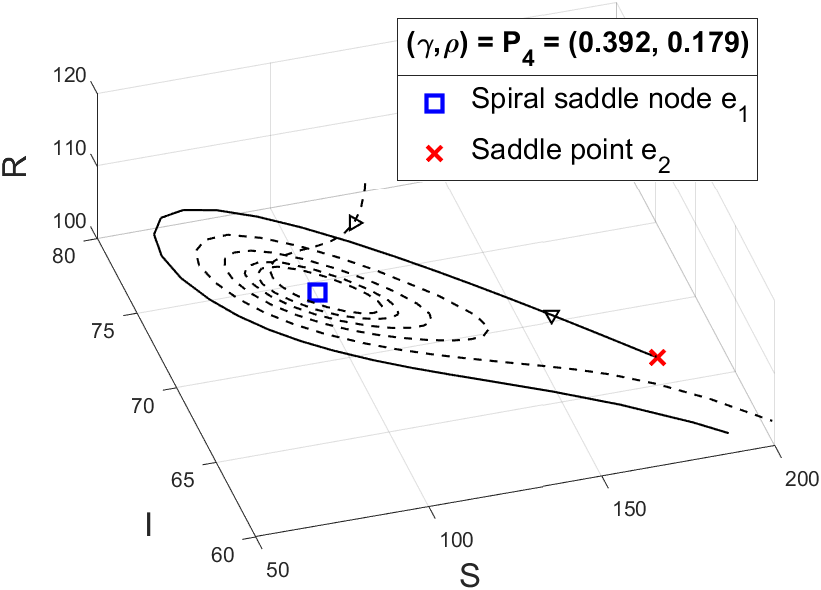}&
		\includegraphics[width=0.3\textwidth]{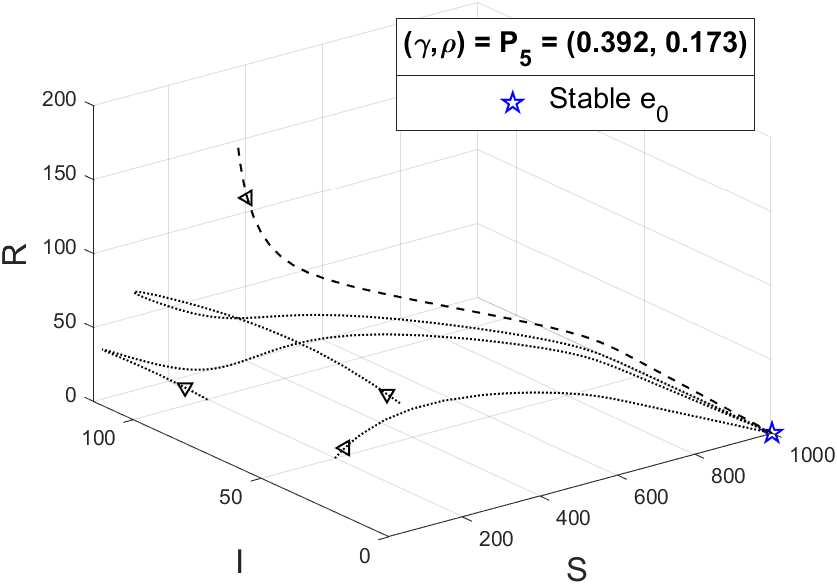}\\	
		{\scriptsize \textbf{(g)} $(\gamma,\rho)=\PP_4$} & {\scriptsize \textbf{(h)} $(\gamma,\rho)=\PP_4$ (magnified)} & {\scriptsize \textbf{(i)} $(\gamma,\rho)=\PP_5$}\\\hline

		\includegraphics[width=0.3\textwidth]{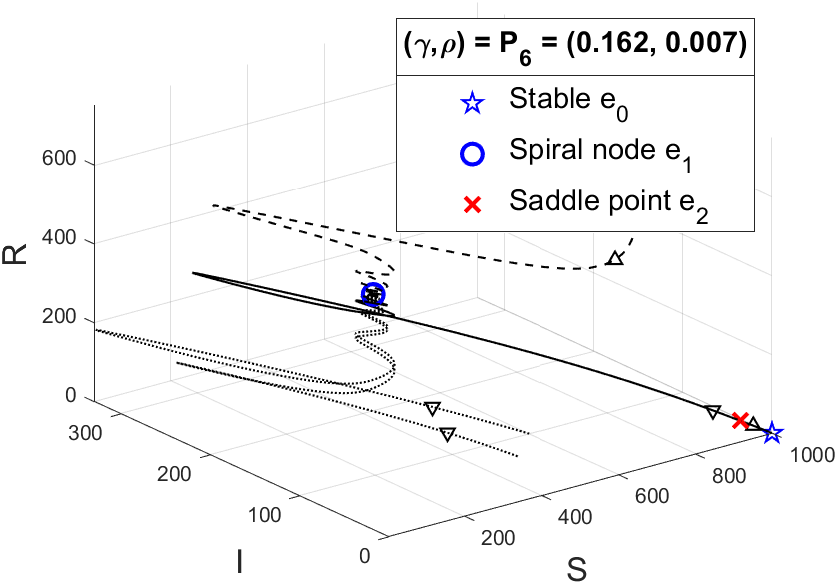}&
		\includegraphics[width=0.3\textwidth]{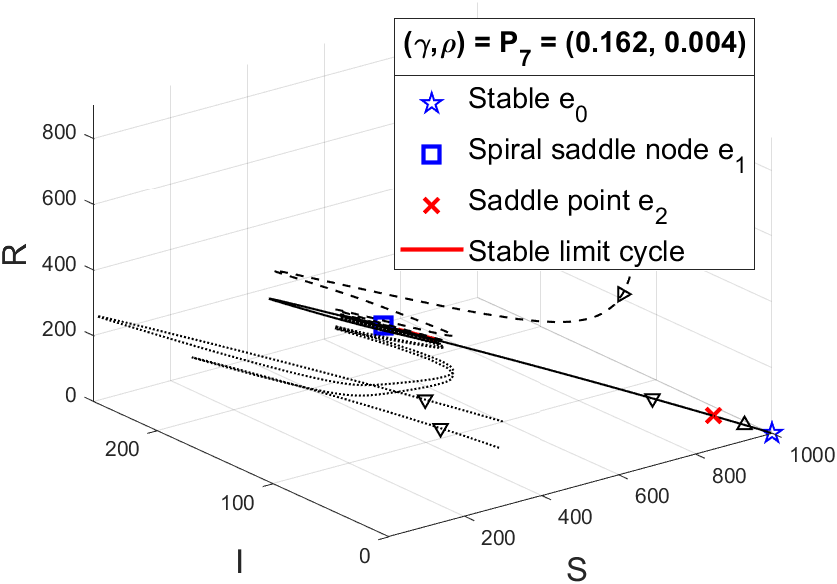}&
		\includegraphics[width=0.3\textwidth]{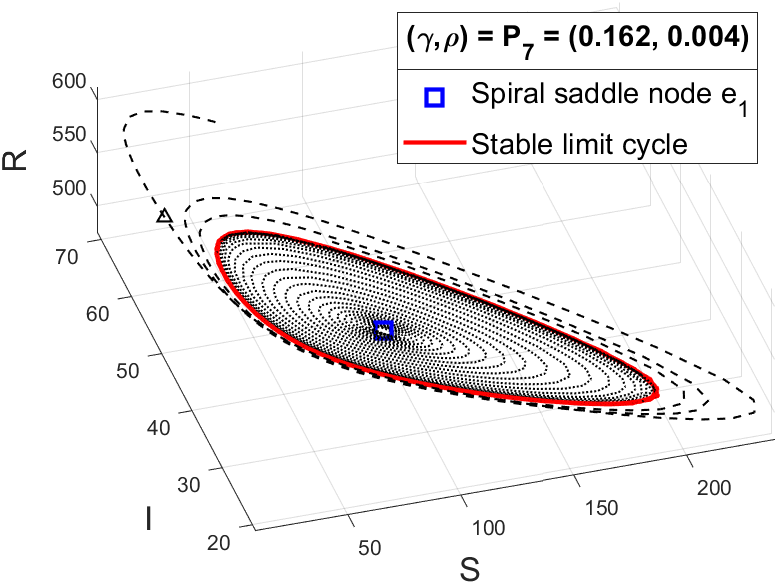}\\	
		{\scriptsize \textbf{(j)} $(\gamma,\rho)=\PP_6$} & {\scriptsize \textbf{(k)} $(\gamma,\rho)=\PP_7$} & {\scriptsize \textbf{(l)} $(\gamma,\rho)=\PP_7$ (magnified)}\\\hline

		\includegraphics[width=0.3\textwidth]{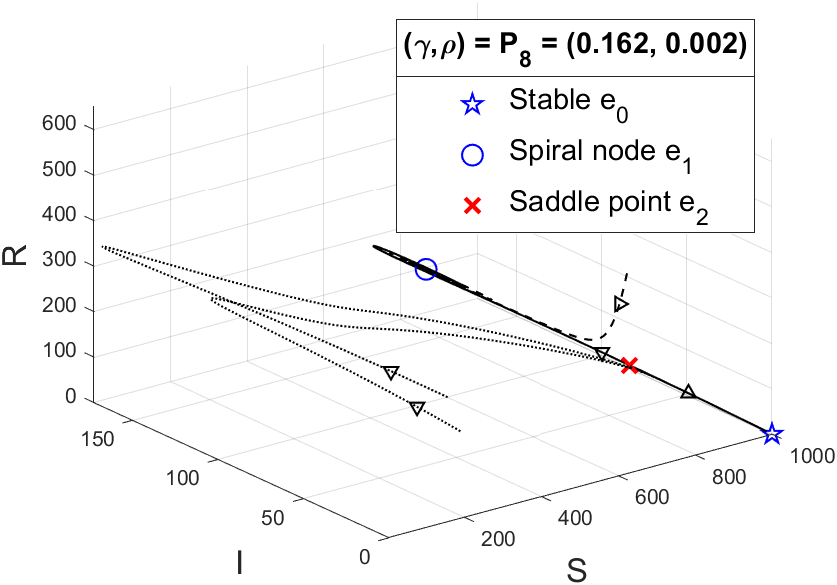}&

\includegraphics[width=0.3\textwidth]{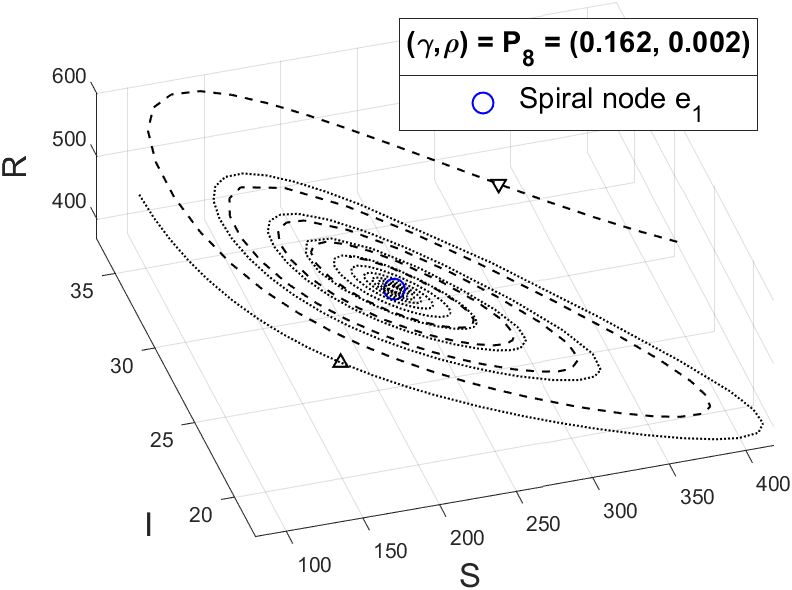}&

\includegraphics[width=0.3\textwidth]{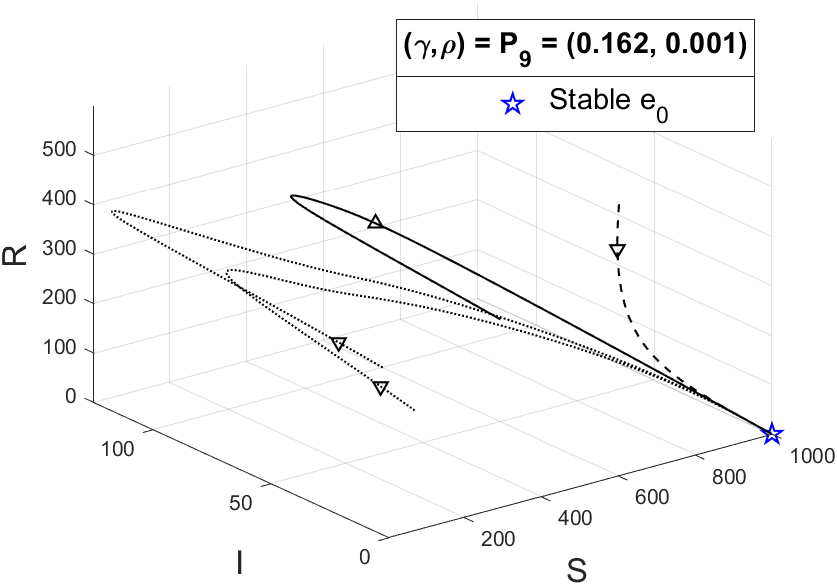}\\	
		{\scriptsize \textbf{(m)} $(\gamma,\rho) = \PP_8$} & {\scriptsize \textbf{(n)} $(\gamma,\rho) = \PP_8$ (magnified)} & {\scriptsize \textbf{(o)} $(\gamma,\rho)=\PP_9$}\\\hline

\end{tabu}}	
		\caption{\label{fig:TE1}Phase portraits of the model \eqref{eq:model} for $\left(\gamma,\rho\right)=\PP_i$, where $i\in\{1,\ldots,9\}$.}
\end{figure}
\begin{figure}\centering
{\tabulinesep=0.5mm
\begin{tabu}{|c|c|c|}\hline
		\includegraphics[width=0.3\textwidth]{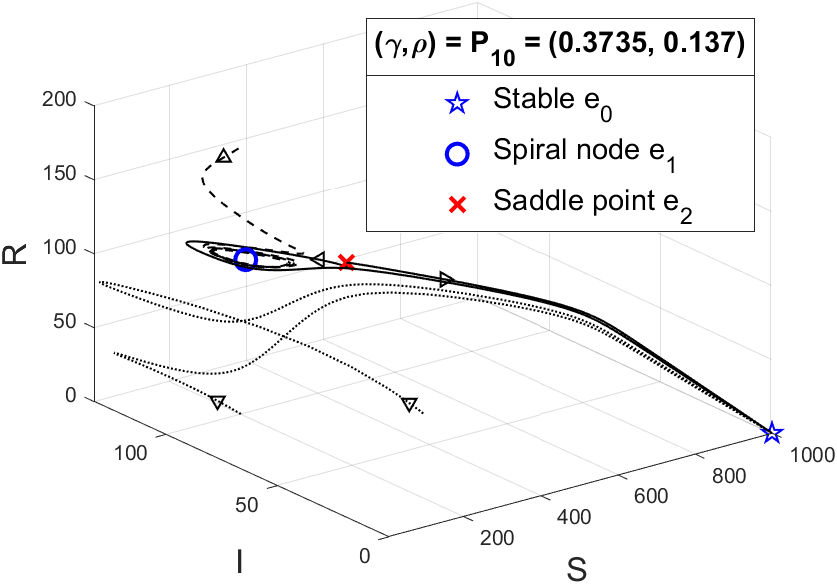}&
		\includegraphics[width=0.3\textwidth]{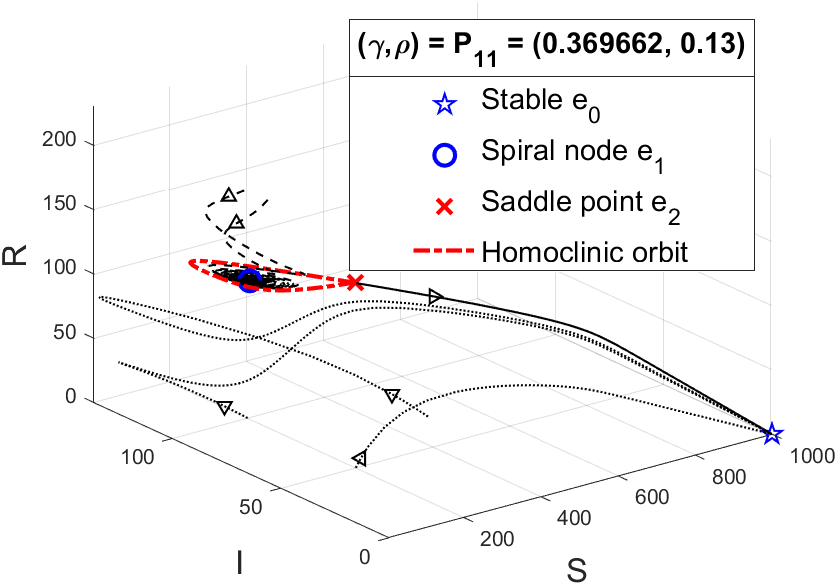}&
		\includegraphics[width=0.3\textwidth]{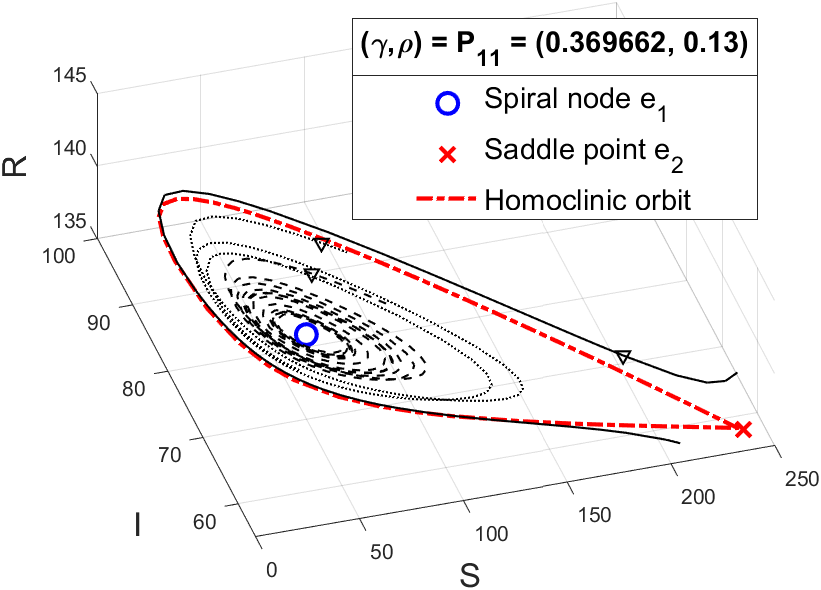}\\
		{\scriptsize \textbf{(a)} $(\gamma,\rho) = \PP_{10}$} & {\scriptsize \textbf{(b)} $(\gamma,\rho) = \PP_{11}$} & {\scriptsize \textbf{(c)} $(\gamma,\rho) = \PP_{11}$ (magnified)}\\\hline
		
		\includegraphics[width=0.3\textwidth]{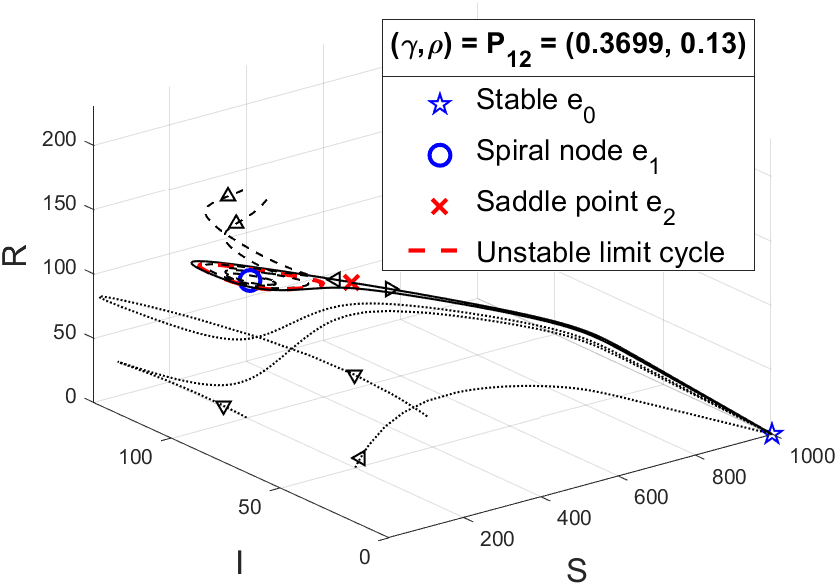}&
		\includegraphics[width=0.3\textwidth]{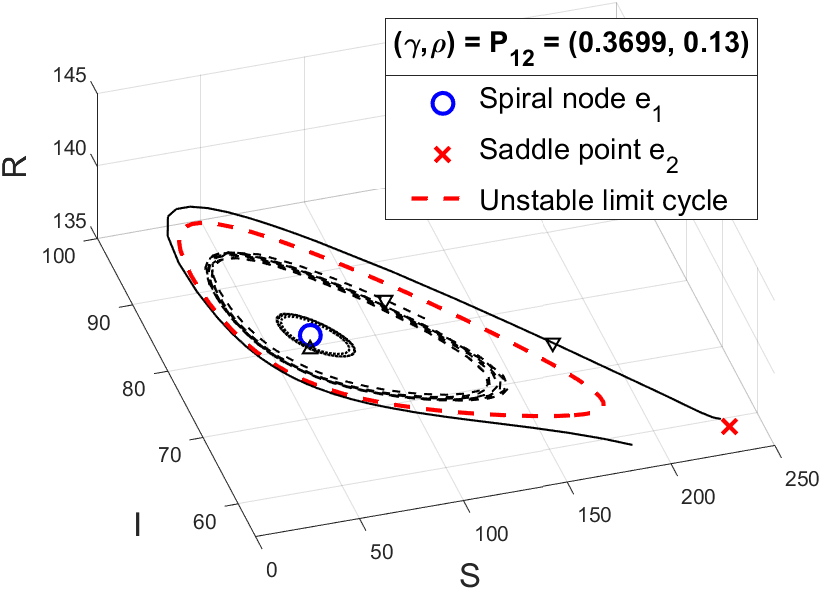}&
		\includegraphics[width=0.3\textwidth]{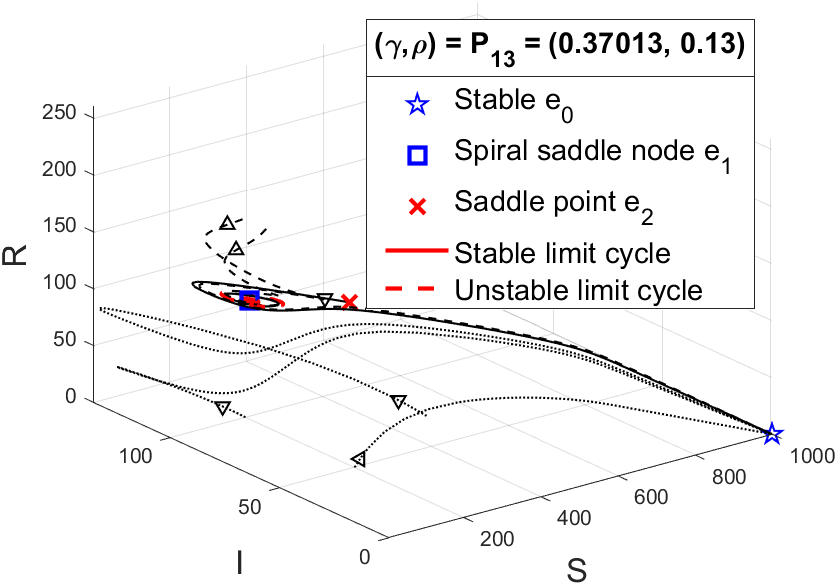}\\
		{\scriptsize \textbf{(d)} $(\gamma,\rho) = \PP_{12}$} & {\scriptsize \textbf{(e)} $(\gamma,\rho) = \PP_{12}$ (magnified)} & {\scriptsize \textbf{(f)} $(\gamma,\rho) = \PP_{13}$}\\\hline

		\includegraphics[width=0.3\textwidth]{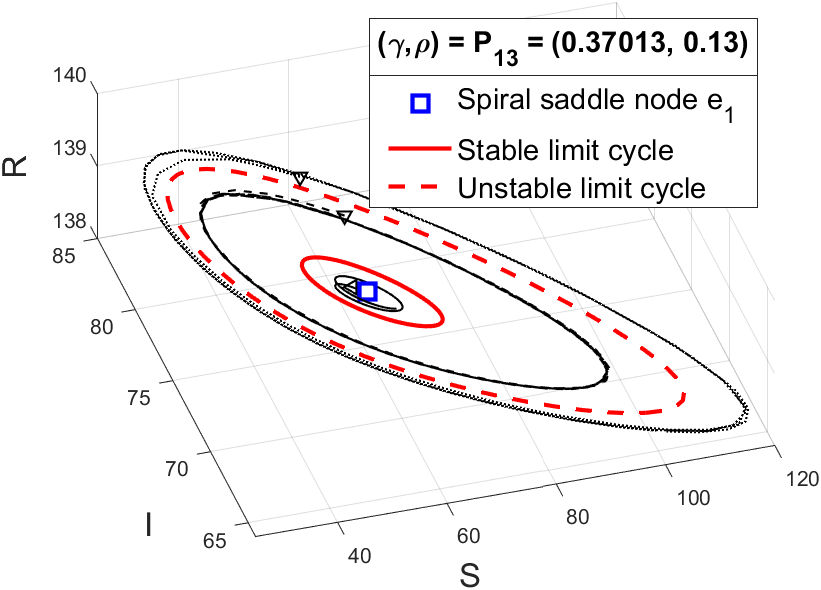}&
		\includegraphics[width=0.3\textwidth]{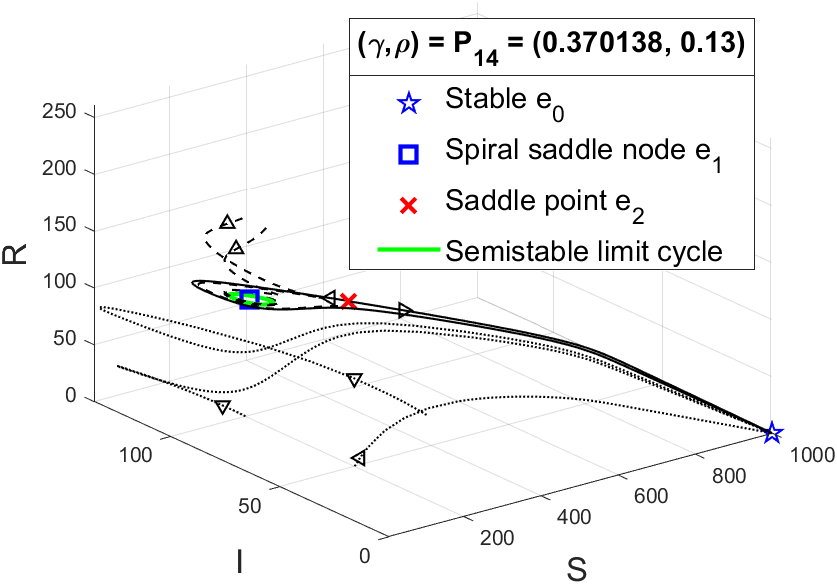}&
		\includegraphics[width=0.3\textwidth]{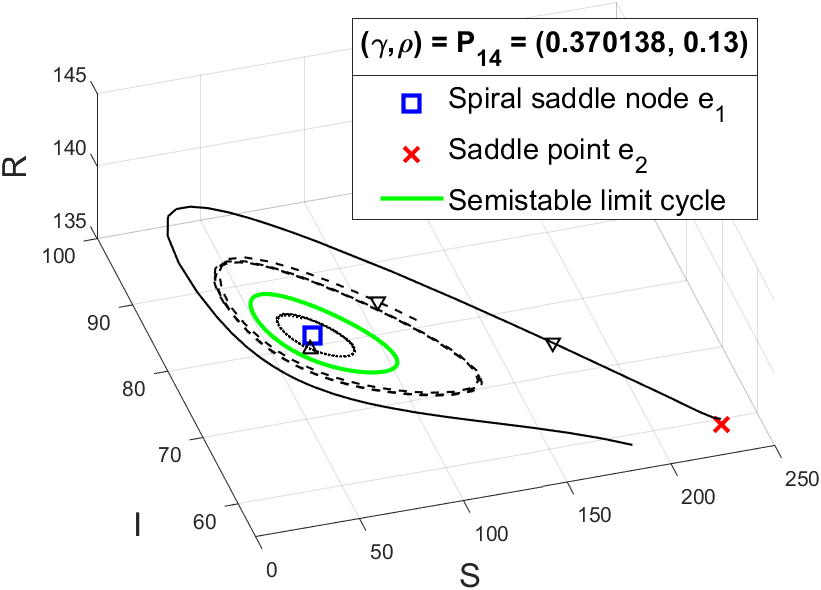}\\	
		{\scriptsize \textbf{(g)} $(\gamma,\rho)=\PP_{13}$ (magnified)} & {\scriptsize \textbf{(h)} $(\gamma,\rho)=\PP_{14}$} & {\scriptsize \textbf{(i)} $(\gamma,\rho)=\PP_{14}$ (magnified)}\\\hline

		\includegraphics[width=0.3\textwidth]{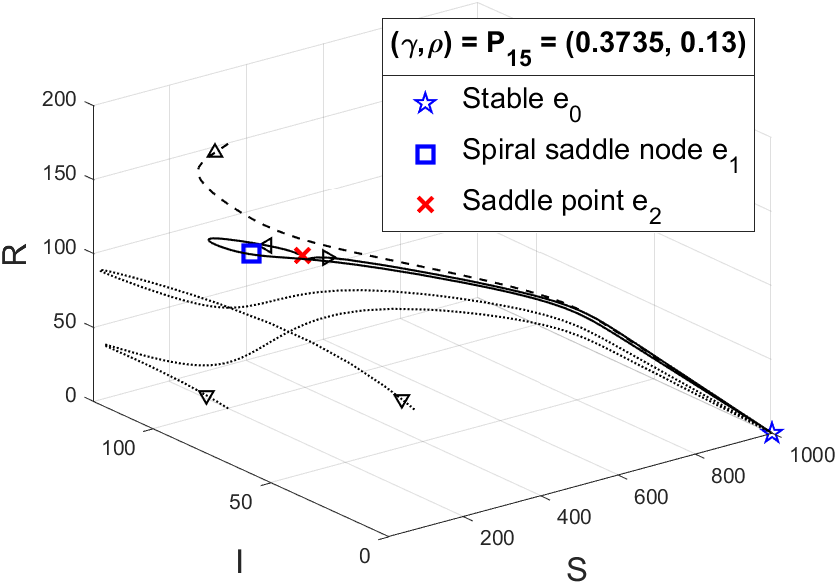}\\	
		{\scriptsize \textbf{(j)} $(\gamma,\rho)=\PP_{15}$}\\\cline{1-1}

\end{tabu}}		
		\caption{\label{fig:TE2}Phase portraits of the model \eqref{eq:model} for $\left(\gamma,\rho\right)=\PP_i$, where $i\in\{9,\ldots,15\}$.}
\end{figure}

In this section, we visualise and describe the model's orbital behaviours at various points $(\gamma,\rho)$ lying in the neighbourhoods of the four bifurcation points, and explain their epidemiological implications. In Figure \ref{fig:BDlocal}, we display magnifications of Figure \ref{fig:BDglobal} in these neighbourhoods. In each neighbourhood, we shall choose a number of specific points $(\gamma,\rho)$ representing a number of qualitatively different orbital behaviours which indicate the occurrence of the respective bifurcation. These behaviours, which we now explain in detail, are all visualised in the model's phase portraits arranged in Figures \ref{fig:TE1} and \ref{fig:TE2}.

\subsection{The dynamical behaviour near $\BT_1$}\label{subsec:BT1}

Figure \ref{fig:BDlocal} (a) shows a magnification of Figure \ref{fig:BDglobal} near the Bogdanov-Takens bifurcation point $\BT_1$, with the addition of a homoclinic bifurcation curve, plotted in red. Let us fix the susceptible individuals' cautiousness level at $\gamma=0.392$, begin with a relatively high value of the hospitals' bed-occupancy rate $\rho$, and describe the topological changes occurring as $\rho$ decreases gradually.
\begin{enumerate}
\item[(1)] At $\rho=0.19$, we have that $(\gamma,\rho)=\PP_1$. Here, the stable disease-free equilibrium $\mathbf{e}_0$ coexist with two endemic equilibria: a stable spiral node $\mathbf{e}_1$ and a saddle point $\mathbf{e}_2$. Therefore, orbits may approach not only the disease-free equilibrium $\mathbf{e}_0$, but also the endemic equilibrium $\mathbf{e}_1$ (Figure \ref{fig:TE1} (a) and (b)). This means that, at this relatively high value of the hospitals' bed-occupancy rate, the disease may persist despite $\cR_0<1$. 
\item[(2)] At $\rho\approx 0.183711$, we have that $(\gamma,\rho)=\PP_2$. At this point, the model undergoes a homoclinic bifurcation. The orbital behaviours remain qualitatively the same, except that a homoclinic orbit now emerges around the saddle endemic equilibrium $\mathbf{e}_2$, acting as a separatrix: orbits inside it approach $\mathbf{e}_1$, others approach $\mathbf{e}_0$ (Figure \ref{fig:TE1} (c) and (d)).
\item[(3)] At $\rho=0.1825$, we have that $(\gamma,\rho)=\PP_3$. Here, the homoclinic orbit has shrunk and become an unstable limit cycle, while remaining a separatrix in the sense previously described (Figure \ref{fig:TE1} (e) and (f)).
\end{enumerate}
As we further decrease $\rho$, we arrive at $\rho=\rho_A\approx 0.181354$, where the stable endemic equilibrium $\mathbf{e}_1$ absorbs the unstable limit cycle while losing its stability, via a subcritical Hopf bifurcation. This leaves no stable endemic equilibrium, and hence the disease's disappearance.
\begin{enumerate}
\item[(4)] At $\rho=0.179$, we have that $(\gamma,\rho)=\PP_4$. Here, no limit cycle exists, and the endemic equilibrium $\mathbf{e}_1$ has become a spiral saddle node. Since no endemic equilibria is stable, orbits approach the disease-free equilibrium $\mathbf{e}_0$ (Figure \ref{fig:TE1} (g) and (h)), meaning that the disease dies out.
\end{enumerate}
Decreasing $\rho$ further, one reaches the backward bifurcation threshold $\rho=\rho_B\approx 0.176117$, where the two endemic equilibria $\mathbf{e}_1$ and $\mathbf{e}_2$ coalesce and disappear via a saddle-node bifurcation, leaving only the stable disease-free equilibrium $\mathbf{e}_0$.
\begin{enumerate}
\item[(5)] At $\rho=0.173$, we have that $(\gamma,\rho)=\PP_5$. Here, no endemic equilibria exist, and orbits still approach the stable disease-free equilibrium $\mathbf{e}_0$ (Figure \ref{fig:TE1} (i)).
\end{enumerate}
From the perspective of the disease's eradication, this analysis highlights the importance of a low bed-occupancy rate. Specifically, for $\gamma=0.392$, in order to guarantee the disease's disappearance, it is necessary to suppress the bed-occupancy rate to below the Hopf bifurcation point $\rho_A$. Notice, however, that $\rho_A$ is larger than the backward bifurcation threshold, i.e., the saddle-node bifurcation point $\rho_B$.

\subsection{The dynamical behaviour near $\BT_2$ and $\GH_2$}

A magnification of Figure \ref{fig:TE1} near the Bogdanov-Takens bifurcation point $\BT_2$ is presented in Figure \ref{fig:BDlocal} (b). Comparing this to Figure \ref{fig:BDlocal} (a), one sees that around $\BT_2$, the model's orbital behaviours are qualitatively the same as those around $\BT_1$.

Let us now turn our attention to the generalised Hopf bifurcation point $\GH_2$, near which a magnification of Figure \ref{fig:TE1} is displayed in Figure \ref{fig:BDlocal} (c). Here let us set $\gamma=0.162$, and again observe the topological changes occurring as $\rho$ is decreased gradually.
\begin{enumerate}
\item[(6)] At $\rho=0.007$, we have that $(\gamma,\rho)=\PP_6$. Here, no limit cycles exist, while two endemic equilibria coexist: the stable spiral node $\mathbf{e}_1$ and the unstable saddle point $\mathbf{e}_2$. Orbits are attracted by both $\mathbf{e}_1$ and the disease-free equilibrium $\mathbf{e}_0$ (Figure \ref{fig:TE1} (j)). Thus, as at $\PP_1$, here we have the possibility of the disease continuing to exist despite $\cR_0<1$.
\end{enumerate}
As $\rho$ is decreased from $0.007$ to $0.004$, it passes through a supercritical Hopf point $\rho=\rho_C\approx 0.002408$, at which $\mathbf{e}_1$ loses stability while ejecting a stable limit cycle.
\begin{enumerate}
\item[(7)] At $\rho=0.004$, we have that $(\gamma,\rho)=\PP_7$. Here, the presence of the stable limit cycle around $\mathbf{e}_1$ implies that the disease could still persist even though the endemic equilibria $\mathbf{e}_1$ and $\mathbf{e}_2$ are both unstable (Figure \ref{fig:TE1} (k) and (l)).
\item[(8)] At $\rho=0.002$, we have that $(\gamma,\rho)=\PP_8$, and we have qualitatively the same behaviours as those at $(\gamma,\rho)=\PP_{6}$ (Figure \ref{fig:TE1} (m) and (n)).
\end{enumerate}
Finally, decreasing $\rho$ further, we arrive at $\rho=\rho_D\approx 0.001573$, where the two endemic equilibria coalesce and disappear in a saddle-node bifurcation.
\begin{enumerate}
\item[(9)] At $\rho=0.001$, we have that $(\gamma,\rho)=\PP_{9}$. Here, the only existing equilibria is the stable disease-free equilibrium $\mathbf{e}_0$ (Figure \ref{fig:TE1} (o)).
\end{enumerate}

In the case of subsection \ref{subsec:BT1}, a complete eradication can already be guaranteed as soon as the bed-occupancy rate becomes lower than the Hopf bifurcation point $\rho_A$. In the present case, by contrast, it is necessary for to suppress the bed-occupancy rate to not merely below the Hopf bifurcation point $\rho_C$, but below the saddle-node bifurcation point, i.e., the backward bifurcation threshold $\rho_D$. Nevertheless, epidemiologically speaking, our conclusion from this analysis is similar, i.e., that a low bed-occupancy rate is necessary for a complete eradication. As strategies to suppress $\rho$, we recommend:
\begin{itemize} 
\item reducing the number of patients having only mild symptoms by optimising self-isolation;
\item transferring near-recovery patients from hospitals to hotels and apartments, so that more beds may be allocated to new patients and queues may be avoided;
\item increasing hospital bed conversions for COVID-19 patients.
\end{itemize}

By keeping the susceptible individuals' cautiousness level $\gamma$ constant, we have demonstrated the importance of having a low bed-occupancy rate $\rho$ for the disease's eradication. In the next subsection, where we describe the orbital behaviours near $\GH_1$, we shall, in turn, fix a specific value of bed-occupancy rate $\rho$ and see the importance of having a high susceptible individuals' cautiousness level $\gamma$.

\subsection{The dynamical behaviour near $\GH_1$}

We now consider the generalised Hopf bifurcation point $\GH_1$, in a neighbourhood of which the bifurcation diagram in Figure \ref{fig:BDglobal} is displayed in Figure \ref{fig:BDlocal} (d), with the addition of a homoclinic bifurcation curve, again plotted in red, and a saddle-node bifurcation of limit cycles curve, plotted in green.
\begin{enumerate}
\item[(10)] At $(\gamma,\rho)=\PP_{10}=(0.3735,0.137)$, the orbital behaviours are qualitatively the same as those at $(\gamma,\rho)=\PP_1$ (Figure \ref{fig:TE2} (a)): no cycles exist, and orbits approach either $\mathbf{e}_1$ or $\mathbf{e}_0$.
\end{enumerate}
Let us now fix $\rho=0.13$, and describe the topological changes occurring as $\gamma$ is increased gradually.
\begin{enumerate}
\item[(11)] At $\gamma=\gamma_A\approx 0.369662$, we have that $(\gamma,\rho)=\PP_{11}$, and that the model undergoes a homoclinic bifurcation: a homoclinic orbit emerges around the saddle endemic equilibrium $\mathbf{e}_2$, being a separatrix: orbits inside it approach the stable endemic equilibrium $\mathbf{e}_1$, others approach the disease-free equilibrium $\mathbf{e}_0$, as at $(\gamma,\rho)=\PP_2$ (Figure \ref{fig:TE2} (b) and (c)). At this low cautiousness level, therefore, we still have the possibility of the disease persisting even though $\cR_0<1$.
\end{enumerate}
Increasing $\gamma$, the homoclinic orbit shrinks and becomes an unstable limit cycle, without abandoning its role as a separatrix.
\begin{enumerate}
\item[(12)] At $\gamma=0.3699$, we have that $(\gamma,\rho)=\PP_{12}$, and that the orbital behaviours are as at $(\gamma,\rho)=\PP_3$ (Figure \ref{fig:TE2} (d) and (e)).
\end{enumerate}
As $\gamma$ is increased from $0.3699$ to $0.37013$, it passes through a subcritical Hopf bifurcation point $\gamma=\gamma_B\approx 0.370127$, where stable endemic equilibrium $\mathbf{e}_1$ loses stability and ejects a stable limit cycle. Here we again have a situation where, although no stable endemic equilibrium exist, the disease could still persist due to the presence of a stable limit cycle.
\begin{enumerate}
\item[(13)] At $\gamma=0.37013$, we have that $(\gamma,\rho)=\PP_{13}$, and that two limit cycles coexist, with opposite stabilities. Orbits near $\mathbf{e}_1$ approach the stable limit cycle, as also those in between the two limit cycles, while orbits outside the unstable limit cycle approach the disease-free equilibrium $\mathbf{e}_0$ (Figure \ref{fig:TE2} (f) and (g)).
\item[(14)] At $\gamma=\gamma_C\approx 0.370138$, we have that $(\gamma,\rho)=\PP_{14}$, and that the two limit cycles coalesce in a saddle-node bifurcation of limit cycles, resulting in a single semistable limit cycle, orbits inside of which approach the limit cycle, while others approach $\mathbf{e}_0$ (Figure \ref{fig:TE2} (h) and (i)).
\item[(15)] At $\gamma=0.3735$, we have that $(\gamma,\rho)=\PP_{15}$, and that the semistable limit cycle no longer exists, so that at $(\gamma,\rho)=\PP_{15}$, orbits are attracted only by the disease-free equilibrium $\mathbf{e}_0$ (Figure \ref{fig:TE2} (j)). It is only at this stage that we are able to guarantee the disease's complete disappearance.
\end{enumerate}
Therefore, for $\rho=0.13$, we have seen that the disease's eradication can only be guaranteed when $\gamma$ exceeds the backward bifurcation threshold $\gamma_C$. As strategies to increase $\gamma$, we recommend:
\begin{itemize}
\item optimising the use of media as tools to educate the public on the risks from COVID-19 and the efforts for prevention;
\item continuing the campaign and enforcement of strict health protocols, so as to help breaking transmission chains.
\end{itemize}

\section{Conclusions and future research}\label{sec:Conclusions}

We have studied a mathematical model for the spread of COVID-19, which incorporates as two main parameters the susceptible individuals' cautiousness level $\gamma$ and the hospitals' bed-occupancy rate $\rho$.  A rectangular region exists on the $\gamma\rho$-plane where $\cR_0<1$, the transcritical bifurcation at $\cR_0=1$ is backward, and four codimension-two bifurcation points exist: two Bogdanov-Takens bifurcation points and two generalised Hopf bifurcation points. Our analysis near each bifurcation point has revealed the complex phenomena through which the model's asymptotic behaviour shifts from endemic to disease-free, which involves the births and disappearances of stable and unstable limit cycles and homoclinic orbits. From an epidemiological viewpoint, the analysis confirms the significance of the two parameters for the eradication of COVID-19. Indeed, the latter can be achieved, provided that susceptible individuals are sufficiently cautious of the disease's spread ---and thus implement the appropriate health protocols--- and that serious efforts are made to keep the hospitals' bed-occupancy rate at a manageable level.

As already noted in \cite{YongOwenHoseana}, the model studied in the present paper is much simplified, and so is modifiable in a number of ways, such as by introducing more compartments and the possibility of reinfection, as realised in \cite{YongHoseanaOwen2}. In addition, since it is quite natural to suspect the nonlinear incidence rate $\beta S I/(1+\gamma S)$ to be a main reason for the emergence of the complex behaviour studied in this paper, one could try replacing it with alternative forms of nonlinear incidence rate \cite{HethcoteDriessche,LiuHethcoteLevin,XiaoRuan}, such as $$\frac{\beta S I}{1+\gamma S^p},\quad \frac{\beta S I}{1+\gamma I^q},\quad \frac{\beta S I}{1+\gamma_1 S^p +\gamma_2 I^q},\quad \beta S^p I^q,$$
and investigate how the dynamical behaviour of the resulting model compares to that of the present model.


\end{document}